\documentclass{article}[12pt]
\hoffset - 2.1cm
\voffset - 1.2cm
\usepackage{amssymb}

\usepackage[english,francais]{babel}

\newtheorem{e-proposition}[theorem]{Proposition}

\newtheorem{e-definition}[theorem]{Definition\rm}

\begin{document}
\thispagestyle{empty}
\centerline{Statistics}
\vskip 3.5cm

 \begin{center}
{\LARGE Nonparametric inference of a trend using functional data}\\

\vskip 4mm

{\large David Degras}
\vskip 1mm

{\it \small Laboratoire de Statistique Th\'eorique et Appliqu\'ee - Universit\'e Paris 6\\
Bo\^ite 158, 175 rue du Chevaleret, 75013 Paris, France.}


\medskip

{\small Received *****; accepted after revision +++++\\
Presented by Paul Deheuvels}
\end{center}
\vskip 7mm
\selectlanguage{english}
\hrule \vskip 3mm

{\small \noindent{\bf Abstract} \vskip 1.5\baselineskip
Let $X=\{ X(t), t \in [0,T] \} $ be a second order random process of which $n$ independent realizations are observed on a fixed grid of $p$ time points. Under mild  regularity assumptions on the sample paths of  $X$, we show the asymptotic normality of suitable nonparametric estimators of the trend function $\mu = \mathbb{E}X$ in the space $C([0,T])$ as $n,p \to\infty$ and, using Gaussian process theory,  we derive  approximate simultaneous confidence bands for $\mu$. 

\vskip 0.5\baselineskip

\noindent{\bf R\'esum\'e} \vskip 0.5\baselineskip \noindent
{\bf Inf\'erence  non param\'etrique d'une tendance avec des donn\'ees fonctionnelles.}  
Soit $X$=$\{ X(t),  t  \in  [0,T] \}$ un processus al\'eatoire du second ordre dont on observe $n$ r\'ealisations ind\'ependantes sur une grille de $p$ points d\'eterministes. Sous de faibles conditions de r\'egularit\'e sur les trajectoires de $X$, nous prouvons la normalit\'e asymptotique d'estimateurs non param\'etriques de la tendance $\mu = \mathbb{E}X$ dans l'espace $C([0,T])$ lorsque $n,p \to\infty$, 
puis nous obtenons des bandes de confiance simultan\'ees approch\'ees pour $\mu$  \`a l'aide de la 
th\'eorie des processus Gaussiens.}
\vskip 4mm \hrule
\quad
\vskip 15mm

{\noindent \bf 1. \hskip 2mm Introduction}\\
\vskip 3mm

In various application fields such as internet traffic monitoring, medical imagery, or signal processing, 
modern technology has allowed to collect data routinely from population samples
with a high temporal and/or spatial resolution. 
 Indeed, such datasets should be viewed as (collections of) curves or functions 
rather than as high-dimensional vectors; they are thus commonly termed \emph{functional data}. 
(See \cite{FerratyVieu2006,RS2005} for\break a comprehensive introduction to functional data analysis.) 
Typical functional data may be modeled as
\vskip 2mm
\noindent \rule[1.2ex]{1.2cm}{0.1mm}

{\small {\it \noindent Email address:}  {\verb david.degras@upmc.fr } (David Degras). 
\vskip 4mm
{\it \noindent Preprint submitted to the Acad\'emie des sciences \hfill November 23, 2008}}
\pagebreak

\noindent observations of independent realizations 
$X_1, \ldots, X_n $ of a second order random process  $X=\{ X(t), t \in \mathcal{D} \}$ 
at fixed design points $t_1, \ldots ,t_p $, where $ \mathcal{D}$ denotes a continuous temporal and/or spatial domain. 
In this framework, the observed data are 
\begin{equation}\label{model}
Y_{ij} = X_i(t_j) + \varepsilon_{ij}, \quad 1\leq i \leq n, \: 1\leq j \leq p ,
\end{equation}
  where the $ \varepsilon_{ij}$ are mean zero random variables (r.v.) representing potential measurement errors.

The trend function $\mu = \mathbb{E}X$ often appears as a population mean response function, which motivates its inference. 
The nonparametric regression literature contains several results on the asymptotic properties of estimators of $\mu$ as the  sample sizes $n$ and $p$ go to infinity.  
For instance when $ \mathcal{D}=[0,1] $, mean-square convergence rates of kernel and spline estimators  can be found  in \cite{BR2006,Cardot2000,D-Jallet2005,HartWehrly1986}.  
When  $\mathcal{D}$ is a compact metric space, \cite{Degras2008} gives
a universal consistency result as well as the asymptotic normality  
of all usual regression estimators in the sense of finite dimensional distributions and of 
the space $L_2(\mathcal{D})$, with an application to simultaneous confidence intervals. 
The task of building (nonparametric and simultaneous) confidence bands for $\mu$, 
which proves useful in various problems of prediction, model diagnostic, or calibration 
(e.g. \cite{ABH1989,KSSY1984}), has received considerable attention in the classical regression setting (e.g. \cite{ES1993,Xia1998}) but not, to our knowledge, for functional data.

In this Note, we study the model (\ref{model}) in the case where the random process $X$ is indexed by a compact interval $\mathcal{D}=[0, T]$ and has mildly regular sample paths. 
In Section 2, we state the asymptotic normality of suitable nonparametric estimators of $\mu$ 
in the space $C([0,T])$ of all continuous functions on $[0,T]$ as $n,p \to\infty$. 
In Section 3, we use the former results to build approximate simultaneous confidence bands for $\mu$. Finally in Section 4, some potential applications and extensions of our results are discussed.


\vskip 9mm 
{\noindent \bf 2. \hskip 2mm Asymptotic normality of nonparametric estimators}\\

We state here the assumptions made on the random process $X$ and on the model (\ref{model}) of Section 1.
 \begin{enumerate}
\item[(A.1)] $X$ is mean-square continuous on $\mathcal{D}=[0,T]$. \vspace*{-2.6mm}
\item[(A.2)] The sample paths $X(\omega,\cdot)$  
are almost surely (a.s.) variation-bounded, with their total variation bounded by $B(\omega)$, where $B$ is a r.v. with finite variance; or  \vspace*{-2.6mm}
\item[(A.2')] $| X(\omega,s) - X(\omega,t) | \leq C | s-t |^{\beta } $ a.s. for some positive 
constants $C$ and $\beta $.  \vspace*{-2.6mm}
\item[(A.3)] $\mu$ has two bounded derivatives on $[0,T]$.  \vspace*{-2.6mm}
\item[(B.1)]  The data form a triangular array : $Y_{ij}$\,$=$\,$Y_{ij}(n)$, $t_j$\,$=$\,$t_j(n)$, and $p$\,$=$\,$p(n)$, with $p(n)\to\infty$ as $n\to \infty$.  \vspace*{-7mm}
 \item[(B.2)] The random errors $ \varepsilon_{ij}$ are mutually independent and independent 
 of the $X_i$; they have mean zero and common variance $\sigma^2 \geq 0$.  \vspace*{-2.6mm}
\item [(B.3)]  The $t_{j}$ are ordered ($0\leq  t_1< \ldots <t_p \leq T $) and they
have a quasi-uniform repartition, i.e., writing $t_0=0$ and $t_{p+1}=T$, it holds that 
 $\frac{\max_{0 \leq j \leq p}(t_{j+1}-t_j )}{ \min_{1 \leq j \leq p-1}(t_{j+1}-t_j )} =\mathcal{O}(1)$ as $n,p\to\infty $. 
\end{enumerate}
\vskip 2pt
\noindent Note that (A.2') implies (A.2). Also, (B.3) ensures that the $t_j$ are regularly spaced in $[0,T]$; it is fulfilled e.g. when the $t_j$ are equally spaced or are generated by a regular probability density function (p.d.f.).

For each  assumption (A.2) and (A.2'), we now introduce a suitable nonparametric estimator of $\mu$ and give its asymptotic distribution.  
Under (A.2), we use the interpolation-type estimator of \cite{Clark1977}, denoted by $\widehat{\mu}_C$,
 with a boundary correction. 
We  recall here its definition. Let $\overline{Y_j}=(\sum_{i=1}^n Y_{ij})/n$ for $  1 \leq j \leq p,$ and let $\overline{Y}(t)$ be the process obtained by linear interpolation of the $(t_j, \overline{Y_j})$ such that $ \overline{Y}(t)=  \overline{Y_1}$ if $t\leq t_1$ and $\overline{Y}(t)= \overline{Y_p} $ if $t\geq t_p$.
 The estimator $\widehat{\mu}_C$ is the convolution of a kernel function $K$ with $\overline{Y}$ :
\begin{equation}\label{Clark-estimator}
\widehat{\mu}_C(t) = \frac{1}{\int_{ 0}^{ T}  K_h(t-u) \mathrm{d}u} 
\int_{ 0}^{ T} K_h(t-u) \overline{Y}(u) \mathrm{d}u , \qquad K_h(\cdot) = K(\cdot/h)/h .
\end{equation}
\pagebreak

\noindent For convenience we take $K$ as a symmetric, compactly supported, Lipschitz-continuous p.d.f.. The real $h>0$ is a fixed bandwidth. 
Following \cite{Pollard1990}, we say that a sequence $(Z_n)$ of random elements of $C([0,T])$ converges weakly to a limit $Z$ in $ C([0,T])$ if $\mathbb{E}\varphi(Z_n) \to \mathbb{E}\varphi (Z)$ as $n\to\infty$ for all uniformly continuous functional $\varphi$ on  $C([0,T])$ equipped with the sup-norm. 
We denote by $R$ the covariance function of the process $X$ and by $\mathcal{G}(0,C)$ any Gaussian process indexed by $[0,T]$ with mean zero and covariance $C$. 
We are now in position to state the weak convergence of  $\widehat{\mu}_C$  in $C([0,T])$ as $n,p\to\infty$ (recall that $p=p(n)$).  
\vskip 1mm

{\noindent \bf Theorem 2.1} \hskip 2mm 
{\it Assume that (A.1),(A.2),(B.1)--(B.3) hold and that $h=h(n,p)\to 0$ and $ph^2\to\infty$ as $n,p\to\infty$. 
Then $n^{1/2} ( \widehat{\mu}_C - \mathbb{E}\widehat{\mu}_C )$
converges weakly to $ \mathcal{G}(0,R)$  in $C([0,T])$.   
If in addition (A.3) holds and $n=o(p)$,  $nh^2 \to 0  $ as $n,p\to\infty$, then 
$ n^{1/2} ( \widehat{\mu}_C - \mu )$ also converges to $ \mathcal{G}(0,R) $ in $C([0,T])$.  }

\vskip 2mm

Next, we address the case where $X$ satifies (A.2') (H\"older continuity). We consider the local linear estimator, denoted here by $\widehat{\mu}_{L}$ and defined by 
\begin{equation}
\widehat{\mu}_{L}(t) = \widehat{\theta}_0, \qquad ( \widehat{\theta}_0,  \widehat{\theta}_1) = \mathrm {argmin}_{(\theta_0, \theta_1)} \sum_{j=1}^p \Big( \overline{Y_j} -  \theta_0 - (t_j-t) \theta_1   \Big)^2 \: K_h(t_j - t)  .
\end{equation}

{\noindent \bf Theorem 2.2} \hskip 2mm 
{\it Assume that (A.1),(A.2'),(B.1)--(B.3) hold and that $h=h(n,p)\to 0$ and $ph^2\to\infty$ as $n,p\to\infty$. 
Then $n^{1/2} ( \widehat{\mu}_L - \mathbb{E}\widehat{\mu}_L )$
converges weakly to $ \mathcal{G}(0,R)$  in $C([0,T])$.   
If in addition (A.3) holds, $n=o(p^2)$,  and $nh^4 \to 0 $ as $n,p\to\infty$, then 
$ n^{1/2} ( \widehat{\mu}_L - \mu )$ also converges to $ \mathcal{G}(0,R) $ in $C([0,T])$.  }

\vskip 2pt
\noindent \emph{Remarks.}
\begin{enumerate}
\item[(1)] The proofs of the former theorems are similar and rely on the following steps: 
(i) note that the estimator is linear in the data; 
(ii)  use the functional central limit theorem 10.6 of \cite{Pollard1990} 
for the estimator applied to the data without noise $X_i(t_j)$; 
(iii) show that under 
condition $ph^2 \to \infty $, 
the estimator applied to the 
errors $\varepsilon_{ij}$ becomes
negligible in probability before $n^{-1/2}$ as $n,p\to\infty$, 
uniformly over $[0,T]$; 
(iv) impose additional conditions on $\mu$ and on the joint rates of $n,p$ and $h$ 
to make the asymptotic bias of the estimator as  $o(n^{-1/2})$ uniformly over $[0,T]$ 
(in particular the rate $n=o(p)$ used in Theorem 2.1 is only needed to control the boundary effects in the bias of  $\widehat{\mu}_C$).\vspace*{-2.5mm}
  \item[(2)] Under (A.2), we can prove asymptotic normality in $C(0,T])$ only for $\widehat{\mu}_C$. 
This is because $\widehat{\mu}_C$, as opposed to more classical estimators, has the remarkable feature of preserving monotonicity and thus satisfies (A.2) like $X$, which makes the step (ii) of the theorem proof straightforward. 
On the other hand, under the stronger assumption (A.2') the asymptotic normality in $C(0,T])$ can be obtained for various classical kernel or projection estimators, as well as for $\widehat{\mu}_C$. 
The choice of $\widehat{\mu}_L$ here was motivated by the popularity of this estimator 
and by its good bias properties. \vspace*{-2.5mm}
\item[(3)] The condition $ph^2\to \infty$ can be dropped in both theorems if there is no noise in the data ($\sigma=0$). Besides, the results carry over to the case of correlated errors, e.g. of autoregressive or mixing type. \vspace*{-2.5mm}
  \item[(4)] Theorem 2.2 corrects a mistake in the Section 4 of \cite{Yao2007} which gives $(nph)^{1/2}$ as the normalizing rate for the weak convergence of $ \widehat{\mu}_L(t) - \mathbb{E} \widehat{\mu}_L(t) $. (The condition (C$3^\ast$) ($ph\to 0$ as $n \to \infty$)  
  of this paper does not produce a well-defined estimator for large $n$). On the other hand, our normalizing rate $n^{1/2}$ is consistent with the variance rate $n^{-1}$ found in the literature.  

\end{enumerate}

\vskip 9mm 
{\noindent \bf 3. \hskip 2mm Simultaneous confidence bands}\\

We build here approximate simultaneous confidence bands for $\mu$ at the level $1-\gamma \in (0,1)$. 
First assume that either Theorem 2.1 or 2.2 applies, i.e. that     
 $n^{1/2}(\widehat{\mu}-\mu)$ converges weakly in $C([0,T])$, where  $\widehat{\mu}$ denotes the corresponding estimator (2) or (3).
Assume also that the covariance function $R$ is non-degenerate and  
let $\widehat{R}(t,t)$ be any uniformly consistent estimator of $R(t,t)$ with respect to $t \in [0,T]$. 
With Slutsky, one sees that $n^{1/2} ( \widehat{\mu}_{np} - \mu) / \widehat{R}(t,t)^{1/2}$ 
converges to  $Z= \mathcal{G}(0,\rho)$ in $C([0,T])$, where $\rho$ is the correlation function of $X$.
It suffices then to apply a classical result of \cite{LS1970} to get that 
\begin{equation}\label{LShepp1970}
\lim_{\lambda \to \infty } \lambda^{-2} \log \mathbb{P} \bigg\{ \sup_{t\in [0,T]}Z(t) > \lambda \bigg\} = - \bigg( 2 \sup_{t\in [0, T]} \rho(t,t) \bigg)^{-1} = - \frac{1}{2}\: .
\end{equation}
Finally, use  (\ref{LShepp1970}) along with the inequality $ \mathbb{P} \{ \sup_{t\in [0,T]}|Z (t)| > \lambda\} \leq 2 \mathbb{P}\{ \sup_{t\in [0,T]}Z(t) > \lambda\}$  to derive the following approximate simultaneous confidence bands for $\mu$ :    
\begin{equation}
\widehat{\mu}(t) \pm \left(- 2\log (\gamma/2) \, \widehat{R}(t,t)  /n \right)^{1/2} \quad \big( 0\leq t \leq T \big).
\end{equation}

\vskip 9mm 
{\noindent \bf 4. \hskip 2mm Discussion}\\

The asymptotic normality results presented in this paper provide a new tool 
for making simultaneous inference on a trend function $\mu$ 
in the context of functional data. 
They can plausibly be extended to the framework of multivariate and/or vector-valued random processes  
and to the inference of derivatives of $\mu$. 
To the best of our knowledge, the simultaneous confidence band procedure of Section 3 is the only one appearing in the literature  for functional data. Its implementation 
only requires the estimation of the mean and of the variance of $X$ 
(not the whole covariance $R$) 
and some simulations have indicated its good performances in terms of 
empirical coverage probability. 
It would benefit from additional features such as bias correction or data-driven bandwidth selection. 
The asymptotic normality results may also be applied to constructing tests for $\mu$. 
A linearity test based on these results is currently under study.




\end{document}